\documentclass[a4paper,10pt]{article}

\usepackage[utf8]{inputenc}
\usepackage{amsmath, amsfonts, amssymb, amsthm}
\usepackage{mathrsfs}
\usepackage{fullpage}

\usepackage{tikz-cd} %рисование диаграмм

\usepackage{xcolor} %цветной текст
\usepackage{indentfirst} %первая строка --- красная

\newtheorem{theorem}{Theorem}
\newtheorem{fact}{Fact}

\newcommand{\T}{\mathbb{T}}
\newcommand{\TT}{\mathcal{T}}
\newcommand{\pphi}{\varphi}

\title{Absence of local unconditional structure in spaces of smooth functions on two-dimensional torus}

\author{Anton Tselishchev}

\date{}

\begin{document}

\maketitle

	\begin{abstract}
		Consider a finite collection $\{T_1, \ldots, T_J\}$ of differential operators with constant coefficients on $\mathbb{T}^2$ and the space of smooth functions generated by this collection, namely, the space of functions $f$ such that $T_j f \in C(\mathbb{T}^2)$. We prove that under a certain natural condition this space is not isomorphic to a quotient of a $C(S)$-space and does not have a local unconditional structure. This fact generalizes the previously known result that such spaces are not isomorphic to a complemented subspace of $C(S)$.
	\end{abstract}

\section{Introduction}

It is well known and easy to see that the space $C^k(\T)$ of $k$ times continuously differentiable functions on the unit circle is isomorphic to $C(\T)$. Also, it has long been known that in higher dimensions the situation is different --- already for two dimensions the space $C^k(\T^2)$ is not isomorphic to $C(\T^2)$. 

This fact was first announced in $\cite{Grot}$ and later generalized in many directions (see \cite{Henk, KislFactor, Kisl0, KwaPel, Sid, PelSen, KisSid, KisMaks, KisMaks2, Maks}). However, the most general and natural framework was introduced only in the quite recent paper \cite{KMSpap} (see also the preprint \cite{KMSprepr} for the two-dimensional case).

More specifically, suppose we have a collection $\TT=\{T_1, T_2, \ldots, T_J\}$ of differential operators with constant coefficients on the torus $\T^2$. So, each $T_j$ is a linear combination of operators $\partial_1^{\alpha}\partial_2^{\beta}$. We call the number $\alpha+\beta$ the order of such a differential monomial and the order of $T_j$ is the maximal order among all monomials involved in it. We consider the following seminorm on trigonometric polynomials $f$:
$$
\|f\|_{\TT} =\max_{1\leq j\leq J} \|T_j f\|_{C(\T^2)}.
$$ 

Now define the Banach space $C^{\TT}(\T^2)$ by this seminorm (that is, factorize over the null space and consider the completion). For example, whet $\TT$ consists of all differential monomials of order at most $k$, we get the space $C^k(\T^2)$. 

In the papers \cite{KisMaks, KisMaks2} the following statement was proved. Suppose that all differential monomials involved in any of $T_j$ are of order not exceeding $k$. Let us drop the junior part of each $T_j$ (this means that we drop all monomials whose order is strictly smaller than $k$). If among the remaining senior parts there are at least two linearly independent, then $C^\TT (\T^2)$ is not isomorphic to a complemented subspace of $C(S)$. (We denote by $S$ an arbitrary uncountable compact metric space. According to Milutin theorem, all the resulting $C(S)$ spaces are isomorphic.) However, if all senior parts are multiples of one of them, the situation was unclear.

So, in the preprint \cite{KMSprepr} (and in the paper \cite{KMSpap} for arbitrary dimensions) a refinement of this statement was proved. In order to state it, we need the concept of mixed homogeneity.

Fix some \textit{mixed homogeneity pattern}, that is, a line $\Lambda$ that intersects the positive semiaxes. The equation of such a line is $\frac{x}{a}+\frac{y}{b}=1$ where $a$ and $b$ are positive numbers. We call such a line \textit{admissible} if all multiindices $(\alpha, \beta)$ such that $\partial_1^\alpha \partial_2^\beta$ is involved in one of $T_j$ lie below $\Lambda$ or on it. This means that all such multiindices must satisfy the following inequality:
$$
\frac{\alpha}{a}+\frac{\beta}{b}\leq 1.
$$

Now we define the senior part of $T_j$ as the sum of all differential monomials involved in $T_j$ whose multiindices lie on the line $\Lambda$ and the junior part as the sum of all other monomials of $T_j$. The senior part is denoted by $\sigma_j$ and the junior by $\tau_j$.

Suppose that for some choice of $\Lambda$ there are at least two linearly independent among all senior parts $\sigma_j$. Then it was proved in \cite{KMSpap, KMSprepr} that $C^{\TT}(\T^2)$ is not isomorphic to a complemented subspace of a $C(S)$ space.

However, in a less general setting, this is not the best known statement. For example, in \cite{KislFactor} it was proved that $C^{k}(\T^2)$ is not isomorphic to any quotient space of $C(S)$. The following theorem generalizes this statement in the described setting.

\begin{theorem}
	If for the collection $\TT$ there are at least two linearly independent operators among $\sigma_j$ \emph(for some choice of an admissible line $\Lambda$\emph), then $C^{\TT}(\T^2)$ is not isomorphic to any quotient space of $C(S)$.
\end{theorem}

This is the first result of this paper.

Also, we note that there is another generalization of the theorem from \cite{KMSpap, KMSprepr} (again in a less general setting). In \cite{KisSid} it was proved that if all operators in the collection $\TT$ are differential monomials and at least two senior monomials (with respect to some pattern) are linearly independent, then the space $C^{\TT}(\T^2)$ does not have local unconditional structure.

Following \cite{GorLew} we give the definition. A Banach space $X$ is said to have local unconditional structure if there exists a constant $C>0$ such that for any finite-dimensional subspace $F\subset X$ there exists a Banach space $E$ with $1$-unconditional basis and two linear operators $R: F \rightarrow E$ and $S: E\rightarrow X$ such that $SRx = x$ for all $x\in F$ and $\|S\|\cdot \|R\| \leq C$. A basis $\{e_n\}$ is $1$-unconditional if for any numbers $\varepsilon_n$ with $|\varepsilon_n| \leq 1$ and any finitary sequence $(\alpha_n)$ the following inequality holds: $\|\sum \varepsilon_n  \alpha_n x_n| \leq \|\sum \alpha_n x_n\|$.

It is worth noting that $X$ has local unconditional structure if and only if its conjugate $X^*$ is a direct factor of a Banach lattice (see \cite{Pietsch}). Since the space $C(S)$ does have local unconditional structure, the non-isomorphism of $C^{\TT}(\T^2)$ to a complemented space of $C(S)$ would also follow once it is proved that $C^{\TT}(\T^2)$ does not have local unconditional structure. This is exactly the statement of the next theorem.

\begin{theorem}
	If for a collection $\TT$ there are at least two linearly independent operators among $\sigma_j$ \emph(for some choice of an admissible line $\Lambda$\emph), then $C^{\TT}(\T^2)$ does not have local unconditional structure.
\end{theorem}

The main ingredients of our proofs are the same as in \cite{KMSpap, KMSprepr}. We use the new embedding theorem established there together with some facts about $p$-summing operators.

At first, we introduce some definitions. A distribution $f$ on the torus $\T^2$ is called proper if $\hat{f}(s,t)=0$ whenever $s=0$ or $t=0$. Next, we need a notion of Sobolev spaces with nonintegral smoothness:
$$
W_2^{\alpha, \beta} (\T^2) = \{f\in C^\infty(\T^2)': \{(1+m^2)^{\alpha/2} (1+n^2)^{\beta/2}\hat{f}(m,n)\} \in \ell^2 (\mathbb{Z}^2) \}.
$$
Of course, the norm of $f$ in $W_2^{\alpha, \beta}(\T^2)$ is defined as $\|\{(1+m^2)^{\alpha/2} (1+n^2)^{\beta/2}\hat{f}(m,n)\}\|_{\ell^2}$.

Now we state the embedding theorem (see Theorem 0.2 and Remark 1.6 in \cite{KMSpap}) which we are going to use.

\begin{fact}
	Suppose that proper distributions $\phi_1, \ldots, \phi_N$ satisfy the following system of equations\emph:
	\begin{align}
	\label{ET}
	-\partial_1^k\pphi_1 = \mu_0; \qquad \partial_2^l \pphi_j - \partial_1^k \pphi_{j+1}=\mu_j, \quad j=1,\ldots, N-1; \qquad \partial_2^l \phi_N=\mu_N,
	\end{align}
	where $\mu_0, \ldots, \mu_N$ are functions in $L^1(\T^2)$ \emph(or measures\emph). Then
	$$
	\sum_{j=1}^N \|\pphi_j\|_{W_2^{\frac{k-1}{2}, \frac{l-1}{2}} (\T^2)} \lesssim \sum_{j=0}^N \|\mu_j\|.
	$$
\end{fact}

Here (and everywhere in this paper) the symbol $A\lesssim B$ means that there exists some constant $C>0$ such that $A\leq CB$.

Several remarks are in order. First, Theorems 1 and 2 hold also for the torus of arbitrary dimension, $\T^n$. But this fact cannot be derived from $2$-dimensional statements (or at least it is unclear how to do this, see \cite{KMSpap} for some explanations). The proofs in higher dimensions are somewhat similar, however, they are much more technically sophisticated (and even require a different embedding theorem, again, see \cite{KMSpap} and Theorem 1.1 there). So, in this article we restrict ourselves to the two-dimensional case.

In this paper we present the proofs of Theorem 1 and Theorem 2. We start with the first theorem because its proof is easier and contains less technical details (however, the reader will see that the proofs of both theorems are quite similar and similar to the proof from the preprint \cite{KMSprepr}).

We note that in the paper \cite{KisSid} it was also proved (again, in case when all operators in $\TT$ are differential monomials and there are at least two linearly independent operators among their senior parts) that if $C^\TT (\T^2)^*$ is isomorphic to a subspace of a space $Y$ with local unconditional structure, then $Y$ contains the spaces $\ell_\infty^k$ uniformly (again, for the definition see \cite{KisSid}). The same statement can also be proved in our situation, but we do not present the details here, because our main goal is to show that, using the embedding theorem from \cite{KMSprepr}, we can adapt various techniques to a more general context. And although this statement implies Theorem 1, we choose to sacrifice the generality for the sake of simplicity and transparency of presentation.
 
Also, a few words should be said about the notation. As it has already been mentioned, we write $A \lesssim B$ if $A\leq CB$ for some constant $C>0$. It will always be clear from the context from which parameters $C$ can depend and from which it cannot. Besides that, the notation $A\asymp B$ means that $A\lesssim B$ and $B \lesssim A$.

The author is kindly grateful to his scientific advisor, S. V. Kislyakov, for posing these problems, for very helpful discussions during the process of their solution and for great help with editing this text.

\section{Nonisomorphism to a quotient of a $C(S)$-space}

Like it was done in \cite{KMSprepr}, we start our proof of Theorem 1 with some simple but helpful observations.

\subsection{Several reductions}

We denote the space of proper functions in $C^{\TT}(\T^2)$ by $C_0^{\TT}(\T^2)$. It is clear that this space is complemented in $C^{\TT}(\T^2)$ (a projection is given by convolution with some measure), so we can prove Theorem 1 for $C_0^{\TT}(\T^2)$ instead of $C^{\TT}(\T^2)$.

Next, suppose that the admissible line $\Lambda$ is given by the equation $x/a+y/b=1$. Let us show that without loss of generality we may assume that $a$ and $b$ are positive integers. Indeed, according to the conditions of Theorem 1, there are at least two points $(r_1, r_2)$ and $(\rho_1, \rho_2)$ with nonnegative integral coordinates on $\Lambda$. We may assume that $r_1>\rho_1$ and $r_2<\rho_2$. Then the equation of $\Lambda$ can be written in the following form:
$$
\frac{x}{r_1-\rho_1}+\frac{y}{\rho_2-r_2}=\frac{\rho_1}{r_1-\rho_1}+\frac{\rho_2}{\rho_2-r_2}.
$$
Now note that we can shift the line $\Lambda$ (and the whole construction) by a vector with integral coordinates. This means that we can change the collection $\TT$ by the collection $\{T_1\partial_1^u\partial_2^v,\ldots T_J\partial_1^u\partial_2^v\}$. The corresponding spaces 
$$
C_0^{\{T_1, \ldots, T_J\}}(\T^2) \quad \hbox{and} \quad C_0^{\{T_1\partial_1^u\partial_2^v,\ldots T_J\partial_1^u\partial_2^v\}}(\T^2)
$$ 
are isomorphic --- isomorphism is given by the map $f\mapsto \partial_1^u\partial_2^v f$. So, by doing this shift we may assume that the equation of $\Lambda$ is the following:
$$
\frac{x}{r_1-\rho_1}+\frac{y}{\rho_2-r_2}=\frac{\rho_1+u}{r_1-\rho_1}+\frac{\rho_2+v}{\rho_2-r_2}.
$$
If we write this equation in the form $x/a_1+y/b_1=1$, then $a_1$ and $b_1$ are the following:
$$
\rho_1+u+ (\rho_2+v)\frac{r_1-\rho_1}{\rho_2-r_2}\quad \hbox{and} \quad \rho_2+v+ (\rho_1+u)\frac{\rho_2-r_2}{r_1-\rho_1}.
$$
Clearly, we can find positive integers $u$ and $v$ so that these two expressions become integers.

So, we assume that the equation of $\Lambda$ is $x/a+y/b=1$ where $a$ and $b$ are positive integers. We denote their greatest common divisor by $N$ and then all points on $\Lambda$ are of the form $(jm, (N-j)n)$ with $0\leq j\leq N$ (here $m=a/N$ and $n=b/N$ so $m$ and $n$ are coprime).

\subsection{Main construction}

Suppose that $C_0^{\TT}(\T^2)$ is isomorphic to a quotient space of $C(S)$. Denote by $P$ the quotient map, $P:C(S)\rightarrow C_0^{\TT}(\T^2)$. 

Due to the reductions we have done, the senior part of every operator from $\TT$ has the following form:
$$
\sigma_s=\sum_{j=0}^N a_{sj}\partial_1^{jm}\partial_2^{(N-j)n}.
$$
We note that the space $C^{\TT}(\T^2)$ depends only on the linear span of operators in $\TT$ so we can change our collection if these changes do not affect its linear span.

Now we consider the matrix $(a_{sj})$. Suppose $j_0$ is the smallest index such that $a_{sj_0}\neq 0$ for at least one $s$. Without loss of generality we may assume that $a_{1j_0}\neq 0$. Then, multiplying $T_1$ by a constant and subtracting a multiple of $T_1$ from other operators, we can ensure that $a_{1j_0}=-1$ and $a_{sj_0}=0$ for every $s>1$. By the assumption of the theorem, there exists $j_1$ such that $a_{sj_1}\neq 0$ for some $s>1$. Again, without loss of generality we assume that $a_{2j_1}=1$ and $a_{sj_1}=0$ for all $s>2$.

Therefore, we have two operators, $T_1$ and $T_2$, whose senior parts are linearly independent. Then for simplicity we denote the coefficients of their senior parts by $a_j$ and $b_j$ respectively, that is
\begin{align*}
\sigma_1=\sum_{j=0}^N a_j\partial_1^{jm}\partial_2^{(N-j)n};\\
\sigma_2=\sum_{j=0}^N b_j\partial_1^{jm}\partial_2^{(N-j)n}.
\end{align*}
Moreover, $T_1$ is the only operator in $\TT$ that involves the differential monomial $\partial_1^{j_0m}\partial_2^{(N-j_0)m}$ and $T_2$ is the only operator in $\TT$ besides maybe $T_1$ that includes the monomial $\partial_1^{j_1m}\partial_2^{(N-j_1)m}$.

Consider the embedding of the space $C_0^{\TT}(\T^2)$ in $C_0^{T_1, T_2}(\T^2)$ (denote it by $i$). Next, we can embed this space into $W_1^{T_1, T_2}(\T^2)$. We denote this embedding by $g$. Here, clearly, the spaces $C_0^{T_1, T_2}(\T^2)$ and $W_1^{T_1, T_2}(\T^2)$ are defined by the seminorms $\max\{\|T_1f\|_{C(\T^2)}, \|T_2 f\|_{C(\T^2)}\}$ and $\max\{\|T_1f\|_{L^1(\T^2)}, \|T_2 f\|_{L^1(\T^2)}\}$ respectively and consist only of proper functions. We note that operator $g$ is $1$-summing, this follows easily from the Pietsch factorization theorem. A good reference on the theory of $p$-summing operators is the book \cite{Wojt} (see Chapter III.F there).

Next, we are going to construct an operator $s$ from $W_1^{T_1, T_2}(\T^2)$ into $W_2^{\frac{m-1}{2},\frac{n-1}{2}}(\T^2)$. Again, the construction will be very similar to that in \cite{KMSprepr} with certain simplifications.

First, we need the following simple fact.

\begin{fact}
	The system \emph{(\ref{ET})} with proper measures \emph(or $L^1$ functions\emph) $\mu_j$ is solvable if and only if the following relation holds true\emph:
	\begin{align}
	\label{solv}
	\sum_{j=0}^N \partial_1^{jk}\partial_2^{(N-j)l}\mu_j = 0.
	\end{align}
\end{fact}
The proof is quite easily done by induction and can be found in \cite{KMSprepr} (see Lemma 2.1 there).

Now take any $f\in W_1^{T_1, T_2}(\T^2)$ and consider the pair of functions $(f_1, f_2)=(T_1f, T_2f)$. Clearly, they satisfy the equation $T_2 f_1-T_1 f_2=0$. This is a differential equation and now we rewrite it in a different form. In order to do this, we note that if $\alpha/a +\beta/b < 1$, then we can express the differential monomial $\partial_1^\alpha \partial_2^\beta$ in terms of $\partial_1^a$ and $\partial_2^b$, using Fourier multipliers:
$$
\partial_1^\alpha \partial_2^\beta f=I_{\alpha\beta}\partial_1^a f + J_{\alpha\beta}\partial_2^b f,
$$
where $I_{\alpha\beta}$ and $J_{\alpha\beta}$ are Fourier multipliers with the following symbols:
$$
\frac{(iu)^{\alpha+a}(iv)^{\beta}}{(iu)^{2a}\pm (iv)^{2b}} \qquad \hbox{and} \qquad \pm\frac{(iu)^{\alpha}(iv)^{\beta+b}}{(iu)^{2a} \pm (iv)^{2b}},
$$
respectively. By this we mean that they act on a function $g\in L^1_0(\T^2)$ by multiplying its Fourier coefficients $\hat{g}(u,v)$ by these expressions. The choice of a sign $\pm$ is determined by the condition $(-1)^a=\pm (-1)^b$, so that the denominators do not vanish when $u$ and $v$ are not equal to zero. In \cite{KMSprepr} it was proved that such multipliers are bounded on $L^1_0(\T^2)$.

\begin{fact}
	The Fourier multipliers $I_{\alpha\beta}$ and $J_{\alpha\beta}$ defined as above are bounded on $L^1_0(\T^2)$.
\end{fact}

Using these multipliers, we can write the junior parts of operators $T_1$ and $T_2$ in the following form:
$$
\sum_{\alpha, \beta} c_{\alpha\beta} (I_{\alpha\beta}\partial_1^a + J_{\alpha\beta}\partial_2^b).
$$
Therefore, we can regroup the terms in the expression $T_2 f_1 - T_1 f_2$ and rewrite it as
$$
\sum_{j=0}^N \partial_1^{jm}\partial_2^{(N-j)n} \mu_j = 0,
$$
where the $\mu_j$ are precisely the functions $b_j f_1-a_j f_2$ when $j\neq 0, N$, $\mu_0$ is equal to $b_0 f_1-a_0f_2$ plus some linear combination of the operators $J_{\alpha\beta}$ applied to $f_1$ and $f_2$, and $\mu_N$ equals $b_N f_1-a_Nf_2$ plus some linear combination of the operators $I_{\alpha\beta}$ applied to $f_1$ and $f_2$.

Now we use Fact 2 and find a solution of the following system of differential equations:

\begin{equation}
\label{sys1}
-\partial_1^m \pphi_1=\mu_0; \qquad \partial_2^n\pphi_j-\partial_1^m \pphi_{j+1}=\mu_j, \quad j=1,\ldots, N-1; \qquad \partial_2^n \pphi_N =\mu_N.
\end{equation}

By Fact 1, all functions $\pphi_j$ lie in $W_2^{\frac{m-1}{2},\frac{n-1}{2}}(\T^2)$. We take the function $\pphi_{j_0+1}\in W_2^{\frac{m-1}{2},\frac{n-1}{2}}(\T^2)$ (it depends linearly on the initial function $f$) and therefore we get a bounded linear operator $s$ from $W_1^{T_1, T_2}(\T^2)$ into $W_2^{\frac{m-1}{2},\frac{n-1}{2}}(\T^2)$. Summing up, we have the following diagram:

$$
C(S)\xrightarrow{P} C_0^\TT(\T^2) \xrightarrow{i} C_0^{T_1, T_2}(\T^2) \xrightarrow{g} W_1^{T_1, T_2}(\T^2) \xrightarrow{s} W_2^{\frac{m-1}{2},\frac{n-1}{2}}(\T^2).
$$

\subsection{Contradiction}

Now we pass to the final part of the proof. We will construct an operator from a finite-dimensional subspace of $W_2^{\frac{m-1}{2},\frac{n-1}{2}}(\T^2)$ to $C(S)$ and use some standard facts from Banach space theory  (mainly, about absolutely summing operators) to get a contradiction. Now let us pass to the details.

Consider the function $v_{pq}:=z_1^p z_2^q \in C^\TT (\T^2)$. We are going to assume that natural numbers $p$ and $q$ satisfy the inequality 
$$
\frac{\delta}{2} q^n \leq p^m \leq \delta q^n,
$$
 where $\delta$ is a small fixed constant (depending of course on our collection $\TT$ but not on $p$ and $q$) which will be chosen later. Also, we will consider only large values of $p$: $p>C$ for some big constant $C$. We always assume that the numbers $p$ and $q$ satisfy these conditions and do not emphasize this later in the present section.

First of all, we note that $\|v_{pq}\|_{C^\TT (\T^2)}\asymp p^{mN}$.

Indeed, if we take any differential monomial $\partial_1^\alpha \partial_2^\beta$ involved in a junior pat of any operator from $\TT$, then we have: $\partial_1^\alpha \partial_2^\beta z_1^p z_2^q = (ip)^\alpha (iq)^\beta z_1^p z_2^q.$ 
 Since this monomial is in a junior part of some operator, the following inequality holds: $\frac{\alpha}{Nm}+\frac{\beta}{Nn}<1$. Therefore, if $\alpha=\alpha_0 m$, then $\beta=(N-\alpha_0-c)n$ for some $c>0$. Hence, the norm of $\partial_1^\alpha \partial_2^\beta v_{pq}$ in $C(\T^2)$ is equal to $p^\alpha q^\beta = p^{\alpha_0 m} q^{(N-\alpha_0-c)n} \asymp p^{m(N-c)}$. Clearly, this quantity can be made arbitrarily smaller than $p^{mN}$ if we make $p$ sufficiently large.

On the other hand, if we apply any differential monomial involved in the senior part of one of the operators (which is of the form $\partial_1^{jm} \partial_2^{(N-j)n}$) to $v_{pq}$, we get a function whose norm is $p^{jm} q^{(N-j)n} \asymp p^{mN}$. Moreover, if $j>j_0$, then $p^{jm} q^{(N-j)n} \asymp \delta^j q^{nN}$ and this quantity can be made arbitrarily smaller than $p^{j_0m} q^{(N-j_0)n} \asymp \delta^{j_0} q^{nN}$ if we make $\delta$ small, therefore $\|T_1 v_{pq}\|_{C(\T^2)} \asymp p^{mN}$. All these facts easily imply that indeed $\|v_{pq}\|_{C^\TT (\T^2)} \asymp p^{mN}$.

Similarly, $\|v_{pq}\|_{C^{T_1, T_2}(\T^2)} \asymp p^{mN}$ and $\|v_{pq}\|_{W_1^{T_1, T_2}(\T^2)} \asymp p^{mN}$. Therefore, we consider functions 
$$
w_{pq}:=\frac{v_{pq}}{p^{mN}}.
$$

By the discussion above, we have $\|w_{pq}\|_{C_0^\TT (\T^2)} \asymp 1$. Hence, there exist functions $f_{pq} \in C(S)$ such that $P(f_{pq})=w_{pq}$ and $\|f_{pq}\|_{C(S)} \leq C$. Besides that, we see that $T_1 w_{pq} = c_{pq} v_{pq}$ and $T_2 w_{pq} = d_{pq} v_{pq}$ where $|c_{pq}|, |d_{pq}| \asymp 1$.

Next, we need to solve the system of differential equations (\ref{sys1}). We recall that $\mu_0$ equals $b_0 c_{pq} v_{pq} - a_0 d_{pq} v_{pq}$ plus some linear combination of the operators $I_{\alpha\beta}$ and $J_{\alpha\beta}$ applied to $T_1 w_{pq}$ and $T_2 w_{pq}$. Therefore, we write it in the following form:
$$
\mu_0 = \xi_{pq} c_{pq} v_{pq} + \eta_{pq} d_{pq} v_{pq} + (b_0 c_{pq} v_{pq} - a_0 d_{pq} v_{pq}).
$$

It is easy to see that $\xi_{pq}, \eta_{pq} = O(p^{-\varepsilon})$ for some small fixed $\varepsilon>0$. Indeed, we simply recall that the symbol of any Fourier multiplier $I_{\alpha\beta}$ is of the form 
$$
\frac{(ip)^{\alpha+a} (iq^\beta)}{(ip)^{2a}\pm (iq)^{2b}}.
$$
The absolute value of this expression can be estimated by 
$$
\Big| \frac{(ip)^{\alpha+Nm} (iq^\beta)}{(ip)^{2Nm}\pm (iq)^{2Nn}} \Big| \asymp \frac{p^\alpha q^\beta}{p^{Nm}} \asymp \frac{p^\alpha \cdot p^{\frac{m}{n} \beta}}{p^{Nm}}.
$$
Here $\alpha/m + \beta/n < N$, therefore this expression indeed equals $O(p^{-\varepsilon})$. The same is true for all operators $I_{\alpha\beta}$ and $J_{\alpha\beta}$.

Now we find a solution of the system of differential equations (\ref{sys1}). Specifically, we are interested in the function $\pphi_{j_0 + 1}$ (that is how we defined the operator $s$).

If $j_0 = 0$, then we need only the first differential equation to find $\pphi_1$. By construction, $j_0=0$ means that $a_0=-1$ and $b_0 = 0$. Therefore, clearly we have $\pphi_1 = k_{pq} \frac{v_{pq}}{p^m}$ where $|k_{pq}| \asymp 1$.

If $j_0 > 0$, then again by construction $a_0 = b_0 = 0$ and we use the first equation from system (\ref{sys1}) to conclude that 
$$
\pphi_1 = \xi_{pq}^{(1)} \cdot \frac{v_{pq}}{p^m}, \quad \hbox{where} \quad \xi_{pq}^{(1)} = O(p^{-\varepsilon}).
$$
 Note that in this case $|\partial_2^n \pphi_1| = |\xi_{pq}^{(1)}\frac{q^n}{p^m} v_{pq}| \asymp |\xi_{pq}^{(1)} v_{pq}|$. Now, if $j_0 = 1$, then we use the second equation to conclude that $\pphi_2 = k_{pq} \frac{v_{pq}}{p^m}$ with $|k_{pq}| \asymp 1$ (again, in this case $\mu_1 = b_1 c_{pq} v_{pq} - a_1 d_{pq} v_{pq}$ and since $j_0=-1$, we see that $a_1 = -1, b_1 = 0$). If $j_0 > 1$, then we conclude from the second equation that $\pphi_2 = \xi_{pq}^{(2)} v_{pq}$ with $|\xi_{pq}^{(2)}|=O(p^{-\varepsilon})$, etc.
 
 Anyway, we see that the following relation holds for a function $\pphi_{j_0+1}$:
 $$
 \pphi_{j_0+1} = k_{pq} \frac{v_{pq}}{p^m}, \quad \hbox{where} \quad |k_{pq}| \asymp 1.
 $$
 
 Now we emphasize the dependence of $\pphi_{j_0+1}$ on $p$ and $q$ so we denote $\pphi^{(p,q)}:=\pphi_{j_0}$. We see that $\{\pphi^{(p,q)}\}$ is an orthogonal system in $W_2^{\frac{m-1}{2}, \frac{n-1}{2}}$ and
 $$
 \|\pphi^{(p,q)}\|_{W_2^{\frac{m-1}{2}, \frac{n-1}{2}}} \asymp p^{-m} p^{\frac{m-1}{2}}  q^{\frac{n-1}{2}} \asymp p^{-1/2} q^{-1/2}.
 $$
 
 Finally, we consider the finite-dimensional operator $A: W_2^{\frac{m-1}{2}, \frac{n-1}{2}} \to C(S)$ which takes $p^{1/2}q^{1/2} \pphi^{(p,q)}$ to $\alpha_{pq} f_{pq}$ where $(\alpha_{pq})$ is an arbitrary sequence of numbers such that $\sum |\alpha_{pq}|^2 = 1$. Here we take $p$ and $q$ satisfying the previous conditions and such that $p\leq M$ for some big number $M$. To be more precise, the operator $A$ is the composition of the orthogonal projection onto $\mathrm{span}\{\pphi_{p,q}\}_{p<M}$ and the operator we have described.
 
 For any function $g\in W_2^{\frac{m-1}{2}, \frac{n-1}{2}}$ with 
 $$
 g = \sum_{p<M} \varkappa_{pq} p^{1/2} q^{1/2} \pphi^{(p,q)}
 $$
 we have:
 $$
 Ag = \sum_{p<M} \varkappa_{pq} \alpha_{pq} f_{pq}.
 $$
 Then the norm of $A$ can be estimated in the following way:
 $$
 \|Ag\|_{C(S)}\lesssim \sum_{p<M} |\varkappa_{pq}\alpha_{pq}|\leq \Big(\sum |\varkappa_{pq}|^2\Big)^{1/2} \Big(\sum |\alpha_{pq}|^2\Big)^{1/2} \leq \Big(\sum |\varkappa_{pq}|^2\Big)^{1/2}\lesssim \|g\|_{W_2^{\frac{m-1}{2}, \frac{n-1}{2}}},
 $$
 and we conclude that $\|A\| \lesssim 1$.
 
 Besides that, we recall that an operator $g$ (see the diagram in the end of subsection 2.2) is 1-summing and therefore $AsgiP$ is also a 1-summing operator. Since this operator acts on a $C(S)$-space, it is also $1$-integral and therefore its restriction to a finite-dimensional subspace (namely, $\mathrm{span}_{p<M}\{f_{pq}\}$) is 1-nuclear (which simply means that it has a finite trace and its trace can be estimated by its norm), see for example \cite[pp. 218--219]{Wojt}. We have: $\mathrm{tr}(AsgiP)\lesssim \|AsgiP\| \lesssim 1$.
 
 Now we are going to prove that this is false. Recall that by our constructions the operator $sgiP$ takes $f_{pq}$ to $\psi^{(p,q)}$, and $A$ takes $\psi^{(p,q)}$ to $p^{-1/2}q^{-1/2}\alpha_{pq} f_{pq}$. Hence, the operator $AsgiP$ has diagonal form (in the basis $\{f_{pq}\}$):
 $$
 AsgiP(f_{pq}) = p^{-1/2}q^{-1/2} \alpha_{pq} f_{pq}.
 $$
 Since its trace is bounded, we infer that 
 $$
 \Big| \sum p^{-1/2} q^{-1/2}\alpha_{pq} \Big| \lesssim 1.
 $$
 This inequality holds for an arbitrary sequence $(\alpha_{pq})$ with $\sum |\alpha_{pq}|^2=1$, which implies that
 $$
 \sum p^{-1} q^{-1} \lesssim 1.
 $$
 But this is clearly false. Indeed, recall that the number of admissible numbers $q$ here is about $p^{m/n}$ and for every such $q$ we have $q\asymp p^{m/n}$. Therefore, our sum can be estimated from below simply by the following: 
 $$
 \sum_{C<p<M} p^{-1}.
 $$
 Since $\sum p^{-1}$ is a divergent series, this is a contradiction and the proof of Theorem 1 is finished.

\section{Absence of local unconditional structure}

Now we start the proof of another main statement of this paper, Theorem 2. Mainly, this proof involves methods from \cite{KisSid} and also the embedding theorem from \cite{KMSpap}, Fact 1. Precisely like in the proof of Theorem 1, we suppose that the space $C^\TT(\T^2)$ has local unconditional structure.

\subsection{Main constructions}

First, we note that we can do the same reductions as in the previous section where we proved Theorem 1. We consider the space $C^\TT_0 (\T^2)$ instead of $C^\TT(\T^2)$, since the passage to a complemented subspace preserves local unconditional structure. Besides that, we are going to assume that all the additional assumptions from Subsection 2.1 are fulfilled. Next, we define operators $i$, $g$, $s$ in the same way as in Subsection 2.2. We still denote the senior part of $T_j$ by $\sigma_j$ and the junior part by $\tau_j$.

Let us denote by $H$ the collection of differential operators corresponding to all points with integral coordinates on the line $\Lambda$. Then we can consider the embedding $j: C_0^H(\T^2)\to C_0^\TT(\T^2)$ that is a continuous operator (see \cite[Theorem 9.5]{BesIlNik}). So, we have the following diagram:
$$
C_0^H(\T^2)\xrightarrow{j}C_0^\TT(\T^2)\xrightarrow{i}C_0^{T_1, T_2}(\T^2)\xrightarrow{g} W_1^{T_1, T_2}(\T^2)\xrightarrow{s} W_2^{\frac{m-1}{2}, \frac{n-1}{2}}(\T^2).
$$

Recall that the operator $g$ is 1-summing. Now we are going to use the following fact (see \cite{GorLew} or \cite[Sec. 23]{Pietsch}):
\begin{fact}
	Let $X$ be a Banach space having local unconditional structure. Then every $1$-summing operator $T$ from $X$ to an arbitrary Banach space $Y$ can be factored through the space $L^1$, i.e., there is a measure $\mu$ and operators $V:X\to L^1(\mu)$ and $U: L^1(\mu)\to Y^{**}$ such that $UV = \kappa T$, where $\kappa:Y\to Y^{**}$ is the canonical embedding. 
\end{fact}

Using this fact, we get the following commutative diagram:

\begin{tikzcd}
	C_0^H(\T^2) \arrow[r, "j"] & C_0^\TT(\T^2) \arrow[r, "i"] \arrow[rd, "V"] & C_0^{T_1, T_2}(\T^2) \arrow[r, "g"]   & W_1^{T_1, T_2}(\T^2) \arrow[r, "s"] & W_2^{\frac{m-1}{2}, \frac{n-1}{2}}(\T^2) \\
	&                        & L^1(\mu) \arrow[rru, "U"] &             &  
\end{tikzcd}
%https://tikzcd.yichuanshen.de/

Now we can consider the dual diagram:

\begin{tikzcd}
C_0^H(\T^2)^* & C_0^\TT(\T^2)^* \arrow[l, "j^*"'] & C_0^{T_1, T_2}(\T^2)^* \arrow[l, "i^*"']  & W_1^{T_1, T_2}(\T^2)^* \arrow[l, "g^*"'] & W_2^{\frac{m-1}{2}, \frac{n-1}{2}}(\T^2) \arrow[l, "s^*"'] \arrow[lld, "U^*"'] \\
&             & L^\infty(\mu) \arrow[lu, "V^*"'] &             &                        
\end{tikzcd}

The next step of the proof is to construct a specific operator which would take elements of $C_0^H(\T^2)^*$ to elements of the space $W_{1/2}^H (\T^2)$ (which is quasi-Banach; the definition will be given below). This construction is the same as in the paper \cite{KisSid} but for the sake of completeness we repeat it here.

Consider the space $W_2^H(\T^2)$ that is defined by means of the following seminorm (and contains only proper functions):

$$
\|f\|_{W_2^H(\T^2)}=\max\limits_{T\in H} \|Tf\|_{L^2(\T^2)}.
$$

Clearly, this is a Hilbert space. Recall that $H=\{\partial_1^{jm}\partial_2^{(N-j)n}\}_{j=0}^N$. The space $W_2^H(\T^2)$ can be identified with a subspace of $L^2(\T^2)\oplus \ldots \oplus L^2(\T^2)$ (there is $N+1$ copy of $L^2(\T^2)$ here); this identification is given by the map 
$$
f\mapsto (\partial_1^{jm}\partial_2^{(N-j)n} f)_{j=0}^N.
$$

Therefore, we can consider the orthogonal projection from the direct sum of $N+1$ copies of $L^2(\T^2)$ to $W_2^H(\T^2)$. We denote this projection by $P$. We need some properties of these operators, so now we state these properties here. All of them are listed in \cite{KisSid}.

First, it is easy to see how $P$ acts on a natural basis of $L^2(\T^2)\oplus \ldots \oplus L^2(\T^2)$. Suppose that $k=(k_1, k_2)$ is a pair of integers and denote by $\phi_k^l$ the following element of the space $L^2(\T^2)\oplus \ldots \oplus L^2(\T^2)$:

$$
\phi_k^l=(0, 0, \ldots, 0, z_1^{k_1}z_2^{k_2}, 0, 0 \ldots, 0),
$$
where $z_1^{k_1}z_2^{k_1}$ is at the $l$th position, $0\leq l\leq N$. Then the following statement can be proved by simple calculations.

\begin{fact}
	If either $k_1$ or $k_2$ equals $0$, then $P(\phi_k^l)=0$. Otherwise,
	$$
	P(\phi_k^l)=\bar{\lambda}_l\Big( \sum_{j=0}^N |\lambda_j|^2 \Big)^{-1} (\lambda_0 z_1^{k_1}z_2^{k_2}, \ldots \lambda_N z_1^{k_1}z_2^{k_2}),
	$$
	where $\lambda_j=(ik_1)^{jm} (ik_2)^{(N-j)n}$.
\end{fact} 

Now we need to understand how $P$ acts on the space $C_0^H(\T^2)^*$. The space $C_0^H(\T^2)$ can be identified with a subspace of $C(\T^2)\oplus\ldots\oplus C(\T^2)$ (in the same way as $W_2^H(\T^2)$ is identified with a subspace of $L^2(\T^2)\oplus \ldots \oplus L^2(\T^2)$). Therefore, we have:
$$
C_0^H(\T^2)^*=(\mathcal{M}(\T^2)\oplus\ldots\oplus \mathcal{M}(\T^2))/\mathcal{X},
$$
where $\mathcal{X}$ is the annihilator of $C_0^H(\T^2)$ in $C(\T^2)\oplus\ldots\oplus C(\T^2)$, that is,
$$
\mathcal{X}=\Big\{(\mu_0, \mu_1, \ldots, \mu_N): \sum_{j=0}^N \int \partial_1^{jm}\partial_2^{(N-j)n}g\, d\bar{\mu}_j = 0\  \forall\, g\in C_0^H(\T^2) \Big\}.
$$

At this point, formally speaking, we should consider an operator $\Phi_M$ that is convolution with the $M$th Fej\'er kernel in both variables, and the operators $P_M$ such that
$$
P_M (F) = P(\Phi_M \mu_0, \Phi_M \mu_1, \ldots, \Phi_M \mu_N), \quad F\in C_0^H(\T^2)^*,
$$
where $(\mu_0, \ldots, \mu_N)$ is any representative of a functional $F$. This formula is meaningful because $P$ is an \textit{orthogonal} projection and if $(\nu_0, \ldots, \nu_N)$ lies in $\mathcal{X}$, then $(\Phi_M \nu_0, \ldots, \Phi_M \nu_N)$ lies in $\mathcal{X}\cap (L^2(\T^2)\oplus\ldots\oplus L^2(\T^2))$, that is, in the kernel of the projection $P$. Now we state the following fact from \cite{KisSid}.

\begin{fact}
	The operators $P_M: C_0^H(\T^2)\to W_{1/2}^H(\T^2)$ are uniformly bounded in $M$.
\end{fact}

The definition of the space $W_{1/2}^H(\T^2)$ should now be clear from the context.

The proof (modulo some technical details) follows from the theory of singular integrals (and Fourier multipliers) with mixed homogeneity developed in \cite{FabRiv} (we see from the formula in Fact 5 that the components of $P$ are Fourier multipliers with a certain homogeneity) and it is even true that these operators are uniformly of weak type $(1,1)$. Of course, there are some technical differences, for example, in \cite{FabRiv} everything was done for the space $\mathbb{R}^n$ instead of $\T^n$. On the other hand, these differences can be overcome quite easily, again, some details can be found in \cite{KisSid}.

Since all the estimates are uniform in $M$, we omit the letter $M$ in our notation. Now we have the following commutative diagram:

\begin{tikzcd}[column sep=small]
W_{1/2}^H(\T^2) & C_0^H(\T^2)^* \arrow[l, "P"'] & C_0^\TT(\T^2)^* \arrow[l, "j^*"'] & C_0^{T_1, T_2}(\T^2)^* \arrow[l, "i^*"']  & W_1^{T_1, T_2}(\T^2)^* \arrow[l, "g^*"'] & W_2^{\frac{m-1}{2}, \frac{n-1}{2}}(\T^2) \arrow[l, "s^*"'] \arrow[lld, "U^*"'] \\
&	&             & L^\infty(\mu) \arrow[lu, "V^*"'] &             &                        
\end{tikzcd}

Now we are going to use some facts from the theory of $p$-summing operators. A good reference here is \cite{KislOp}. The space $W_{1/2}^H(\T^2)$ is a quasi-Banach space of cotype 2 and $L^\infty(\mu)$ is a space of type $C(K)$. Therefore, $Pj^*V^*$ is a $2$-summing operator (this is a generalization of Grothendieck's theorem; see \cite{KislOp} for details). Hence, $Pj^*i^*g^*s^*$ is also $2$-summing. In the next subsection we are going to show that this is not the case.

\subsection{Final computations and a contradiction}

As in the proof of Theorem 1, we denote by $v_{pq}$ the function $z_1^p z_2^q$. Again, we consider only sufficiently large values of $p$ and assume that the pair $(p, q)$ in question satisfies the following condition:
$$
\frac{\delta}{2}q^n \leq p^m \leq \delta q^n.
$$

We see that 
$$
\|v_{pq}\|_{W_2^{\frac{m-1}{2}, \frac{n-1}{2}}} \asymp p^{\frac{m-1}{2}} q^{\frac{n-1}{2}} \asymp p^m p^{-1/2}q^{-1/2}.
$$

Now denote by $w_{pq}$ the function $\frac{v_{pq}}{\|v_{pq}\|}$. This is an orthonormal system in the space $W_2^{\frac{m-1}{2}, \frac{n-1}{2}}$ and hence it is weakly 2-summable. Therefore, since $Pj^*i^*g^*s^*$ is a 2-summing operator (which by definitions means that it takes weakly 2-summable sequences to 2-summable sequences), we have:
$$
\sum \|Pj^*i^*g^*s^* w_{pq}\|_{W_{1/2}^H}^2 <\infty.
$$

First, let us realize where the operator $j^*i^*g^*s^*$ takes the function $w_{pq}$. Take any function $v_{\tilde{p}\tilde{q}}=z_1^{\tilde{p}}z_2^{\tilde{q}}\in C_0(\T^2)$; linear combinations of such functions are dense in $C_0(\T^2)$. We write:
\begin{equation}
\label{dualop}
\langle v_{\tilde{p}\tilde{q}}, (j^*i^*g^*s^*) w_{pq} \rangle = \langle (sgij) v_{\tilde{p}\tilde{q}}, w_{pq} \rangle = \langle sv_{\tilde{p}\tilde{q}}, w_{pq} \rangle_{W_2^{\frac{m-1}{2}, \frac{n-1}{2}}}.
\end{equation}
Now we need to recall how $s$ acts on the function $v_{\tilde{p}\tilde{q}}$. We need to solve the system of equations (\ref{sys1}) and by definition all functions $\mu_j$ are multiples of $z_1^{\tilde{p}} z_2^{\tilde{q}}$. Therefore, the solution is also a multiple of $z_1^{\tilde{p}} z_2^{\tilde{q}}$ and so we see that $\langle sv_{\tilde{p}\tilde{q}}, w_{pq} \rangle_{W_2^{\frac{m-1}{2}, \frac{n-1}{2}}}\neq 0$ only if $p=\tilde{p}$ and $q=\tilde{q}$.

So, now we need to determine the function $s(v_{pq})$. Recall that in the previous section (where we proved Theorem 1) we showed that $s$ takes $\frac{v_{pq}}{p^{mN}}$ to $k_{pq}\frac{v_{pq}}{p^m}$ where $|k_{pq}|\asymp 1$. Therefore, we have:
$$
s(v_{pq}) = k_{pq} p^{mN} p^{-m} v_{pq}.
$$
Hence, we have the following identity:
$$
\langle sv_{pq}, w_{pq} \rangle = k_{pq} p^{mN}p^{-m}\|v_{pq}\|_{W_2^{\frac{m-1}{2}, \frac{n-1}{2}}}\asymp k_{pq}p^{mN}p^{-m}p^{-1/2}q^{-1/2} = k_{pq} p^{mN}p^{-1/2}q^{-1/2}.
$$
Therefore, we arrive at the following formula for the right-hand side of (\ref{dualop}):
$$
\langle sv_{\tilde{p}\tilde{q}}, w_{pq} \rangle_{W_2^{\frac{m-1}{2}, \frac{n-1}{2}}}=
\begin{cases}
0, \  (p,q)\neq (\tilde{p}, \tilde{q}),\\
k_{pq} p^{mN}p^{-1/2}q^{-1/2}, \ (p,q) = (\tilde{p}, \tilde{q}).
\end{cases}
$$

Recall that the following element of $C(\T^2)\oplus\ldots\oplus C(\T^2)$ corresponds to $v_{pq}\in C_0^H(\T^2)$:
$$
(\partial_1^{jm}\partial_2^{(N-j)n} v_{pq})_{j=0}^N = ((ip)^{jm}(iq)^{(N-j)n} v_{pq})_{j=0}^N.
$$
So, since $p^{mN}\asymp q^{nN}$, we can take the following representative from the equivalence class corresponding to $(j^*i^*g^*s^*)w_{pq}$:
$$
(l_{pq}p^{-1/2}q^{-1/2} v_{pq}, 0, 0, \ldots, 0), \quad \hbox{where} \quad |l_{pq}|\asymp 1.
$$

Finally, we apply the projection $P$ (using the formula from Fact 5; in our case, $|\lambda_j|=p^{jm}q^{(N-j)n}\asymp p^{mN}$ and hence $\bar{\lambda}_l \lambda_k (\sum |\lambda_j|^2)^{-1} \ \asymp 1$). We have:
$$
\sum \|l_{pq} p^{-1/2}q^{-1/2} v_{pq}\|_{L^{1/2}}^2\asymp \sum p^{-1}q^{-1},
$$
and it was already established that this sum is divergent. Therefore, we get a contradiction and the theorem is proved.

\medskip

A. Tselishchev

\medskip

Saint Petersburg Leonard Euler International Mathematical Institute,  Fontanka 27, St. Petersburg 191023, Russia

\medskip

celis-anton@yandex.ru
\end{document}